\title{
High degree graphs contain large-star factors}
\author{ 
Noga Alon
\thanks{
Tel Aviv University, Tel
Aviv 69978, Israel and IAS, Princeton, NJ, 08540, USA. 
Research supported in part by the Israel Science
Foundation, by a USA-Israel BSF grant,
by NSF grant CCF 0832797
and by the Ambrose Monell Foundation.
Email: nogaa@tau.ac.il}
\and Nicholas Wormald
\thanks{Department of Combinatorics and Optimization, University of
Waterloo, Waterloo ON, Canada. Supported by the Canada Research
Chairs Program and NSERC. Email: nwormald@uwaterloo.ca}
}
\date{}
 \documentclass [12pt]{article}
\usepackage{amsmath}
\newtheorem{theo}{Theorem}[section]
\newtheorem{claim}{Claim}[section]

\newtheorem{coro}[theo]{Corollary}

\newcommand{\qed}{\hspace*{\fill} \rule{7pt}{7pt}}

\begin{document}
\maketitle
\begin{abstract}
We show that any finite simple graph with minimum degree $d$
contains a spanning star forest in which every connected component
is of size at least $\Omega((d/\log d)^{1/3})$. This settles a
problem of J. Kratochvil.
\end{abstract}

\section{Introduction}
All graphs considered 
here are finite and simple. A {\em star} is a 
tree with one vertex, the {\em center}, adjacent to all the 
others, which are {\em leaves}. A {\em star factor} of a
graph $G$ is a spanning forest of $G$ in which every connected
component is a star. It is easy to see that any graph with positive
minimum degree contains a star factor in which every component is a
star with at least one edge. Jan Kratochvil \cite{Kr} 
conjectured that if the
minimum degree is large than one can ensure that all stars are
large. More precisely, he conjectured that there is a function
$g(d)$ that tends to infinity as $d$ tends to infinity, so that
every graph with minimum degree $d$ contains a star factor in which
every star contains at least 
$g(d)$ edges. 
Our main result shows that this is
indeed the case, for a function $g(d)$ that 
grows moderately quickly with $d$, as follows.

\begin{theo}
\label{t11}
There exists an absolute positive constant $c$ so that
every graph with minimum degree $d$ contains a star factor
in which
every star has at least $c d/(\log d)^{1/3} $ edges.
\end{theo}

The motivation for Kratochvil's conjecture arises in the running
time analysis of a recent exact exponential time algorithm for the
so called $L(2,1)$-labeling problem of graphs. See \cite{HKKKL}
for more details.

As preparation for the proof of the main result, we prove the
following simpler statement.

\begin{theo}
\label{t12} 
There exists an absolute positive constant $c'$ such that
every $d$-regular graph contains a star factor in which every star
has at least $c'd/\log d $ edges. This is optimal, up to the
value of the constant $c'$.
\end{theo}

Throughout the paper we make no attempt to optimize the absolute
constants. To simplify the presentation we omit all floor and
ceiling signs whenever these are not crucial. We may and will
assume, whenever this is needed, that the minimum degree $d$
considered is sufficiently large.  
All logarithms are in the natural base, unless otherwise specified.

Our notation is standard. In particular, for a
graph $G=(V,E)$ and a vertex $v \in V$, we let $N_G(v)$ denote the set
of all neighbors of $v$ in the graph $G$, and let $d_G(v)=|N_G(v)|$
denote the degree of $v$ in $G$. 
For $X \subset V$, $N_G(X)=\cup_{x \in X}
N_G(x)$ is the set of all neighbors of the members of $X$.

The rest of this short paper is organized as follows. In Section 2 we
present the simple proof of Theorem \ref{t12}, and in Section 3
the proof of the
main result. 
Section 4 contains some
concluding remarks and open problems.

\section{Regular graphs}

\noindent 
{\bf Proof of Theorem \ref{t12}:}\,
Let $G=(V,E)$ be a $d$-regular graph. 
Put $p=(2 + 2\log d)/d$ and let $C$ be a random set of vertices
obtained by picking each vertex of $G$, randomly and independently,
to be a member of $C$, with probability $p$. We will show 
that, with positive probability, some such set $C$ will be a 
suitable choice for the set of centres of the stars in the 
desired star factor. 

For each vertex $v \in V$,
let $A_v$ be the event that either $v$ has no neighbors in $C$ or
$v$ has more than $3pd$ neighbors in $C$. 
By the standard known estimates
for binomial distributions (c.f., e.g., \cite{AS}, Theorem A.1.12),
the probability of each event $A_v$ is at most 
$(1-p)^{d} +(e^2/27)^{2+2 \log d} <1/ed^2$.
Moreover, each event $A_v$ is mutually independent of all 
events $A_u$ except those that satisfy $N_G(v) \cap N_G(u) \neq
\emptyset$. As there are at most $d(d-1)<d^2$ such vertices $u$
we can apply the Lov\'asz Local Lemma (c.f., e.g., \cite{AS}, 
Corollary 5.1.2)
to conclude that with 
positive probability none of the events $A_v$ holds.
Therefore, there is a choice of a set $C \subset V$ so that
for every vertex $v$, $0 < |N_G(v) \cap C| \leq 3pd= 6+ 6\log d.$

Fix such a set $C$, and let $B$ be 
the bipartite graph whose two classes of vertices
are $C$ and $V\setminus C$, where each $v \in C$ is 
adjacent in $B$ to all vertices
$u \in V\setminus C$ which are its neighbors in $G$. 
By the choice of $C$, for every
vertex $v \in C$, $d_B(v) \geq d-6-6 \log d$ and for every vertex
$u \in V\setminus C$, $0<d_B(u) \leq 6 + 6 \log d$. It thus follows by 
Hall's theorem that one can assign to each vertex $v \in C$ 
a set consisting of 
$$
\frac{d-6-6 \log d}{6+6 \log d} >\frac{d}{7 \log d}
$$ 
of its neighbors
in $V\setminus C$,  where no member of $V\setminus C$ 
is assigned to more than one such $v$. (To see this from the standard 
version of Hall's theorem, split each vertex $v$ in $C$ into 
$(d-6-6 \log d)/(6+6 \log d)$ identical `sub-vertices', each with 
the same neighbours in $V\setminus C$ as $v$, and find a matching 
that hits every sub-vertex using Hall's theorem. Then for each 
$v\in C$,  coalesce the subvertices of $v$ back together to form $v$.)
By assigning each unassigned vertex $u$ of $V\setminus C$ 
arbitrarily to one
of its neighbors in $C$ (note that there always is such a neighbor)
we get the required star factor, in which the centers 
are precisely the members
of $C$ and each star contains more than $ d/(7 \log d )$ edges.

It remains to show that the above estimate is optimal, up to a
constant factor. Note that the centers of any 
star factor form a dominating set
in the graph, and thus if the minimum size of a 
dominating set in a $d$-regular graph on $n$ vertices is at least 
$\Omega(n \log d/d)$, then the star factor must 
contain a component of size
at most $O( d/\log d)$. It is not difficult to check that 
the minimum size of a dominating set
in a random $d$-regular graph on $n$ vertices is $\Theta(n  \log d/d)$
with high probability.
In fact, if $c<1$, the expected number of dominating sets of 
size $k=n(c+o(1))(\log d)/d$ tends to 0 for $d$ sufficiently large.   
We give some details of verifying this claim. One can use the 
standard pairing or configuration model of random $d$-regular graphs, 
in which  there are $n$ buckets with $d$ points in each bucket. 
It is enough to prove the result for the multigraph arising from 
taking a random pairing of the points and regarding the buckets 
as vertices (see e.g.~\cite{W} for details). The expected number 
of dominating sets $S$ of size $k$ with $m$ edges from $S$ to $N(S)$ is 
$$
A:={n \choose k} f(k,m)\bigg(\prod_{i=0}^{m-1} 
(kd-i)\bigg)\frac{M(kd-m)M((n-k)d-m)}{M(nd)}
$$
where $f(k,m)$ is the coefficient of $x^m$ in 
$\big( (1+x)^{d}-1\big)^{n-k}$ 
and $M(r)$ is the number of pairings or perfect matchings of an 
even number $r$  of points, i.e.\ $(r-1)(r-3)\cdots 1$. 
In the above formula, the binomial chooses the $k$ buckets of $S$, 
$f(k,m)$ is the number of ways to choose $m$ points in the 
other $n-k$ buckets such that at least one point comes from each 
bucket (so that $S$ dominates the graph), and the next factor counts 
the ways to pair those points with points in buckets in $S$. 
The other factors in the numerator count the ways to pair up the 
remaining points, first within $S$, and then within the rest of 
the graph. The denominator is the total number of pairings. 
We may use  standard methods to see that 
$f(k,m)\le  \big( (1+x)^d-1\big)^{n-k}x^{-m}$ for all real 
$x>0$ (since $f$'s coefficients are all nonnegative). We set 
$x=k/n$ and take $k=nc(\log d)/d$, which is justified by regarding 
$c$ as a function of $n$ that tends to a limit equal to the value 
$c$ referred to in the claim above. This gives
$$
A\le{n \choose k}\big( (1+k/n)^d-1\big)^{n-k}\bigg(\frac{n}{k}\bigg)^m  
\frac{(kd)!}{ (kd-m)!}\cdot\frac{M(kd-m)M((n-k)d-m)}{M(nd)}.
$$
Considering replacing $m$ by $m+2$ shows this expression is 
maximised when $n^2(kd-m)\approx k^2\big((n-k)d-m\big)$ and 
certainly only when $kd-m\sim n\big((c\log d)^2/d\big)$. 
In fact, $kd-m=O\big(n(\log^2 d)/d\big)$ is sufficient for 
our purposes. Fixing such a value of $m$, using   
$n!= (n/e)^nn^{\theta(1)}$ and $M(r)=\Theta((r/e)^{r/2})$, 
and noting that  
\begin{eqnarray*}
{n \choose k}&=& \exp\big(O(n(\log^2 d)/d)\big),\\
(1+k/n)^d-1 &=& d^c\big(1-d^{-c}+O(d^{-1}\log^2 d)\big),\\ 
\bigg(\frac{n}{k}\bigg)^m &=&  
\bigg(\frac{n}{k}\bigg)^{-(kd-m)}\bigg(\frac{n}{k}\bigg)^{kd}
=\bigg(\frac{n}{k}\bigg)^{kd}\exp\big(O(n(\log^3 d)/d)\big),\\ 
((kd-m)/e)^{kd-m} &=&  n^{kd-m}\exp\big(O(n(\log^3 d)/d)\big),\\
((n-k)d-m)^{((n-k)d-m)/2}&=& (nd-2kd+(kd-m))^{((n-k)d-m)/2}\\
&=& (nd)^{((n-k)d-m)/2}\left(1-\frac{2k}{n} 
+\frac{k -m/d}{2n}\right)^{nd(1-2k/n +(k -m/d)/2n)/2}\\
&=& (nd)^{((n-k)d-m)/2}\exp\big(-dk+O(n(\log^2 d)/d)\big) 
\end{eqnarray*}
we find  everything in the upper bound for $A$ cancels except for  
$$
e^{-kd}d^{c(n-k)+(kd-m)/2}(1-d^{-c})^n \exp 
\big(O(nd^{-1}\log^3 d )\big). 
$$
Since $d^c = e^{kd/n}$ and the power of $d$ is absorbed in 
the error term, this equals
$n^{O(1)}\big(1-d^{-c}+O(d^{-1}\log^3 d )\big)^n$, which, 
if $c<1$, tends to 0 for large $d$ as required.
On the other hand, for  $c>1$ and large $d$, $n c (\log d)/d$  
is an upper bound   on the minimum dominating set size in 
all $d$-regular graphs~\cite[Theorem~2.2]{AS}.

An explicit example can be given as well: 
if $d=(p-1)/2$ with $p$ being a prime,
consider the bipartite graph $H$ 
with two classes of vertices $A_1=A_2=Z_p$ in which
$a_ib_j$ forms an edge iff $(a_i-b_j)$ 
is a quadratic non-residue.  A simple
consequence of Weil's Theorem (see, e.g., \cite{Al}, Section 4)
implies that for every set $S$ of at most, say,
$\frac13 \log_2 p$ 
elements of $Z_p$, there are more than $\sqrt p$ members
$z$ of $Z_p$ so that $(z-s)$ is a quadratic residue for all
$s \in S$. This implies that any 
dominating set of $H$ must contain either more than 
$\frac13 \log_2 p$  vertices of 
$A_2$ or at least $\sqrt p$ vertices 
of $A_1$, and is thus of size bigger than
$\frac13 \log_2 p$. (This can in fact be improved to
$(1-o(1)) \log_2 p$, but  as we are not interested in 
optimizing the absolute constants here and in the rest of the paper,
we omit the proof of this stronger statement).
For degrees $d$ that are 
not of the form $(p-1)/2$ for a prime $p$
one can take any spanning $d$-regular subgraph 
of the graph above with the smallest 
$p$ for which  $(p-1)/2  \geq d$. Such a subgraph exists
by Hall's theorem, and any dominating set in it is also 
dominating in the original 
$(p-1)/2$-regular graph, hence it is of size at 
least $\frac{1}{3}\log_2 p$.
By the known results about 
the distribution of primes this prime $p$ 
is $(2+o(1))d$, and we thus get 
a $d$-regular graph on at most $n=(4+o(1))d$ 
vertices in which every dominating
set is of size greater than $\frac13 \log_2 p > \frac13\log_2 d
\geq \Omega \big(n(\log d)/d\big)$.
This completes the proof.  \qed

\section{The proof of the main result}

In this section we prove Theorem \ref{t11}. 
The idea of the proof is based on that of Theorem~\ref{t12}.
Given a graph $G=(V,E)$ with minimum
degree $d$, we wish to define a dominating set 
$C \subset V$ whose members will form the
centers of the star factor, 
and then to assign  many leaves to
each of them. The trouble is that here we cannot pick the set of centers
randomly, as our graph may contain a large set $R$ of vertices of
degree $d$ whose total number 
of neighbors is much smaller than $|R|$, and then
the number of centers in $R$ is limited. 
This may happen if some or all of the neighbors
of the vertices in $R$ have degrees which are much higher than $d$. 
Thus, for example, if our graph is a 
complete bipartite graph with classes of
vertices $R$ and $U$, with $|R|=n-d$ 
and $|U|=d << n-d$, it is better not to 
choose any centers in $R$. In fact, 
it seems reasonable in the general case to force all vertices
of degree much higher than $d$ to be centers, and indeed this is the 
way the proof starts. However, if we then 
have a vertex all (or almost all) 
of whose neighbors have already been declared to be
centers, then this  vertex cannot be a center itself, 
and will have to be a leaf. 
Similarly, if almost all neighbors of a vertex
are already declared to be leaves, then this vertex will have
to become a center.

The proof thus proceeds by declaring, iteratively, some vertices
to be centers and other vertices to be leaves. At the end, if there are
any vertices left, we choose a small subset of them randomly to be 
additional centers. 
The Local Lemma has to be applied to 
maintain the desired properties that
will enable us to apply Hall's theorem at the end to 
a bipartite graph, defined in a way similar to 
that in the proof of Theorem \ref{t12}. 
An additional complication arises from the fact that 
we have to assign time labels to vertices and use them in the
definition of the bipartite graph.  We proceed with the detailed proof.
\vspace{0.4cm}

\noindent
{\bf Proof of Theorem \ref{t11}:}\, 
Let $G=(V,E)$ be a graph with minimum 
degree $d$. We first modify $G$ 
by 
omitting any edge whose two endpoints are of 
degree strictly greater than $d$,
as long as there is such an edge. We thus may and will assume,
without loss of generality,
that every edge has at least one endpoint of degree exactly $d$. Put 
$h=\frac{1}{10}d^{4/3}/(\log d)^{1/3}$, and 
let $H$ (for {\em High}) denote 
the set  of all vertices of degree at least $h$. 
Since each of their neighbors
is of degree precisely $d$, we can apply Hall's theorem and assign 
a set of 
$h/d$ neighbors 
to each of them, so 
that no vertex is assigned twice.   Let $S'$ denote the set of all
the $|H| \cdot\frac{1}{10}d^{4/3}/(\log d)^{1/3}$ 
assigned vertices, and let
$S$ (for {\em Special}) be a 
random subset of $S'$  
obtained by choosing each member of $S'$ to be in $S$ randomly
and independently with probability $1/2$. 
\begin{claim}
\label{c31}
In the random choice of $S$, with positive probability 
the following conditions hold:

\noindent
(i) For each $v \in H$, $|N_G(v) \cap S| 
> \frac{1}{25}d^{1/3}/(\log d)^{1/3}$.

\noindent
(ii) For each $v \in V\setminus H$, $|N_G(v)\setminus  S| \geq d/3.$
\end{claim}
{\bf Proof:}\, For each vertex $v \in H$ let $A_v$ be the event that
condition (i) is violated for $v$. 
Similarly, for each vertex $v \in V\setminus H$
let $B_v$ be the event that condition (ii) is violated for $v$. By the 
standard known estimates for binomial distributions, the probability of
each event  $A_v$ is $\exp\big(-\Omega(d^{1/3}/(\log d)^{1/3})\big)$ 
and that of
each event $B_v$ is $e^{-\Omega(d)}$.
In addition, each event  $A_v$ is independent  of all other events
except for the events $B_u$ for vertices $u$ that have a neighbor among
the $\frac{1}{10}d^{1/3}/(\log d)^{1/3}$ 
vertices of $S'$ assigned to $v$
(note that there are less than $d^{4/3}$ such vertices $u$).
The same reasoning shows that each event $B_v$ is mutually independent 
of all other events $A_u, B_w$ with the exception of
at most $h d<d^3$ events.
The desired result thus follows from the Local Lemma
(with a lot of room to spare). This completes the proof of the claim.

Fix an $S$ satisfying the assertion  of the claim, and
define $G'=G-S=(V', E')$. 
Note that by the above claim, part (ii),
\begin{equation}
\label{e31}
\mbox{for each }~~ v \in V(G')=V\setminus S, ~
|N_{G'}(v)| \geq d/3.
\end{equation}
Note also that by part (i) of the claim, each vertex $v \in H$ has
a set of at least $\frac{1}{25}d^{1/3}/(\log d)^{1/3}$ vertices from $S$
assigned to it, and can thus serve as a center 
of a star of at least that size.

We now construct two sets of vertices $C,L \subset V\setminus S$. 
The set $C$ will 
consist of vertices that are declared to
be {\em centers}, and will serve
as centers of stars in our final star factor. The set $L$ will
consist of vertices that are declared to be {\em leaves} 
in the final factor.  
Note that the vertices in $S$ will not form part 
of these sets; they will also
be leaves in the final star factor, and their 
associated centers will be the 
vertices in $H$ to which they have been assigned, but since we have
already specified their centers we do not need to 
consider them any more.
Initially, define $C=H$ and $L=\emptyset$. 
We will also need a time label
$t(v)$ which will be defined in the following for each vertex in 
$V\setminus (H \cup S)$; in the beginning set $t=0$.

Put $V'=V(G')=V\setminus S$ and define $D=d^{2/3} (\log d)^{1/3}.$ 
Now apply repeatedly the following two rules to define additional
centers and leaves, and assign them time labels.

\begin{itemize}
\item
{\bf (a)} If there is a vertex $v \in V'\setminus (C \cup L)$ 
such that $|N_{G'}(v)\setminus C| \leq D$,
add $v$ to $L$, increase $t$ by $1$, and define $t(v)=t.$
\item
{\bf (b)} If there is a vertex $v \in V'\setminus (C \cup L)$ 
such that $|N_{G'}(v)\setminus L| \leq d/6$,
add $v$ to $C$,  increase $t$ by $1$, and define $t(v)=t.$
\end{itemize}

 The process continues by repeatedly 
applying rules (a) and (b)  in any order  until
there are no vertices left in 
$V'\setminus (C \cup L)$ that satisfy the conditions
in rule (a) or in rule (b). Let $t_0$ denote the value of the 
time parameter $t$ at this point. Actually, since $L$ and 
$C$ will be disjoint, by~(\ref{e31}) no vertex will satisfy 
the conditions in both rules simultaneously. Let $F$ 
(for {\em Free}) denote the set of all
vertices in $V'\setminus (C \cup L)$ remaining once the 
process terminates.
Define $p= 20 (\log d)/d$, and let
$T$ be a random subset of $F$ obtained
by picking each vertex $v \in F$, 
randomly and independently, to be in $T$
with probability $p$.  Assign the vertices  of 
$F\setminus T$ the time labels
$t_0+1, t_0+2, \ldots ,t_0+|F\setminus T|$ in any order.
Finally, assign the 
vertices of $T$ the time labels 
$t_0+|F\setminus T|+1, t_0+|F\setminus T|+2, \ldots ,t_0+|F|$.

In our final star factor, the vertices $H \cup C \cup T$ will serve as
centers, while the remaining vertices, that is, 
those in $S \cup L \cup (F\setminus  T)$,
will serve as leaves. In order to show 
that it is possible to define large stars
with these centers and leaves, we  need the following. 
\begin{claim}
\label{c32}
With positive probability, every vertex of 
$F$ has at least one neighbor in
$C \cup T$, and no vertex $v \in V'-H$ has more than
$2ph =4 d^{1/3} (\log d)^{2/3}$ neighbors in $T$.
\end{claim}
{\bf Proof:}\, 
For each vertex $v \in F$ that does not have any neighbor in $C$, let
$A_v$ be the event that it has no neighbor in $T$. Note that as
$v \in F$, the definition of rule (b) implies that
$|N_{G'}(v)\setminus L|> d/6$, and as it has no neighbor in $C$,
it has more than $d/6$ neighbors in $F$.  
Therefore, the probability that none
of these neighbors is in $T$ is at most $(1-p)^{d/6} <d^{-3}.$
For each vertex $v \in V'\setminus H$, let $B_v$ be the 
event that $v$ has
more than $2ph$ neighbors in $T$. Since the  degree of $v$ in $G'$
is at most $h$, its number of neighbors in $F$ is certainly at most 
$h$, and hence the standard estimates for binomial distributions imply
that the probability  of each event $B_v$ is 
at most $e^{-\Omega(ph)}$ which is
much smaller than, say, $d^{-3}$. 

Note that each event $A_v$ is mutually independent of all other events
$A_u$ or $B_w$ apart from those corresponding 
to vertices $u$ or $w$ that
have a common neighbor with $v$ in $F$, and the number of such 
vertices $u,w$ is smaller than $hd<d^{7/3}$.  Similarly,  
each of the events $B_v$ s independent of all others but
at most
$hd<d^{7/3}$. The claim thus follows from
the Local Lemma.

Returning to the proof of the theorem, fix a choice of $F$ 
satisfying the
assumptions in the last claim. Let $B$ be the bipartite graph with
classes of vertices $(C \setminus H) \cup T$ and 
$ L \cup (F \setminus T)$, in which each
$v \in (C \setminus H) \cup T$ is adjacent to any of 
its neighbors $u$ that lies
in $L \cup (F \setminus T)$ and satisfies 
$t(u) <t(v)$. Note that, crucially,
prospective centers are connected in $B$ 
only to prospective leaves with smaller time labels.

Our objective is to show, using Hall's theorem, that we can assign
to each vertex $v$ in $(C\setminus H) \cup T$ 
some $\Omega(d^{1/3}/(\log d)^{1/3})$
neighbors of $v$ (in $B$, and hence also 
in $G'$) from $L \cup (F \setminus T)$,
such that each vertex in $L \cup (F \setminus T)$ 
is assigned at most once. 
To do so,
we first establish 
several simple properties of the bipartite graph $B$
that follow from its construction. 
\begin{claim}
\label{c33}
The following properties hold.

\noindent
(i) For each vertex $u \in L$, $|N_B(u) \cap (C\setminus H)|
 \leq D=d^{2/3} (\log d)^{1/3}.$

\noindent
(ii) For each vertex $u \in L$, 
$|N_B(u) \cap T| \leq 4 d^{1/3} (\log d)^{2/3}.$

\noindent
(iii) For each vertex $u \in F\setminus T$,
$|N_B(u)|=|N_B(u) \cap T| \leq 4 d^{1/3} (\log d)^{2/3}.$

\noindent
(iv) For each vertex $v \in C \setminus H$, $d_B(v) \geq d/6.$

\noindent
(v) For each vertex $v \in T$, 
$|N_B(v)| \geq D-4 d^{1/3} (\log d)^{2/3}>  D/2=
\frac{1}{2} d^{2/3} (\log d)^{1/3}.$
\end{claim}
{\bf Proof:}\,

\noindent
(i) By the definition of rule (a), each $u \in L$ can have at most
$D$ neighbors  with time labels exceeding $t(u)$, and therefore
can have at most that many neighbors in $B$.

\noindent
(ii) This follows immediately from the condition in 
Claim \ref{c32} that $F$ was chosen to satisfy.

\noindent
(iii) By the definition of the graph $B$, 
the vertices in $F\setminus T$ are
joined in $B$ only to vertices of $T$, as these are the only 
vertices with bigger time labels. Therefore, 
$|N_B(u)|=|N_B(u) \cap T|$ for each $u \in F\setminus T$, 
and the claimed 
upper estimate
for this cardinality follows from Claim \ref{c32}.

\noindent
(iv) By the definition of rule (b), 
each vertex $v \in C\setminus H$ satisfied
$|N_{G'}(v)\setminus L| \leq d/6$ at the point of being added to $C$. 
Since by (\ref{e31}), $|N_{G'}(v)| \geq d/3$, it follows that
at that time, $v$ had at least $d/3-d/6=d/6$ neighbors in $L$.
As all these leaves have smaller time labels than $v$, it is joined
in $B$ to all of them.

\noindent
(v)  If $v \in T$, then $v \in F$, and thus, by the definition of
rule (a), $|N_{G'}(v) \setminus C| > D$. 
Vertices in $T$ are given the largest
time labels, so in the graph $B$, $v$ is joined
to all members of $N_{G'}(v) \setminus C$ 
except for those that lie in $T$. However, by
the condition in
Claim \ref{c32}, at most $4 d^{1/3} (\log d)^{2/3}$ of these vertices
are members of $T$, implying the desired estimate. This 
completes the proof of the claim.
\begin{coro}
\label{c34}
For each subset $X \subset (C \setminus H) \cup T$,
$|N_B(X)| \geq |X|\cdot \frac{1}{16}d^{1/3}/(\log d)^{1/3}.$
\end{coro}
{\bf Proof:}\, 
If at least half the elements of $X$ belong to $C \setminus H$, then,
by Claim \ref{c33} (iv), 
the total number of edges of $B$ incident with them is at least 
$\frac{1}{2}|X| \cdot \frac{1}{6}d$. By 
the first observation in the proof of part (iii) of the claim,
these edges are not incident in $B$ with any member
of $F \setminus T$. By part (i) of the claim,
at most  $D=d^{2/3} (\log d)^{1/3}$ of these edges 
are incident with any one
vertex in $L$. Thus, in this case,
$|N(X)| \geq \frac{1}{2} |X|\cdot \frac{1}{6}d \cdot 1/D
= |X|\cdot \frac{1}{16}d^{1/3}/(\log d)^{1/3}$, 
providing the required estimate.

Otherwise, at least half of the vertices of $X$ lie in $T$. By Claim
\ref{c33}, part (v), the total number of edges of $B$ incident with 
them is greater than $\frac{1}{2} |X|\cdot \frac{1}{2}D $. 
By parts (ii) and (iii)
of the claim, each  neighbor of these 
vertices in $L \cup (F \setminus T)$
is incident with at most  $ 4 d^{1/3} (\log d)^{2/3}$ of 
these edges, implying that in this case
$$
|N(X)| \geq \frac{|X|}{2} 
\frac{D}{2} \frac{1}{ 4 d^{1/3} (\log d)^{2/3}}
=\frac{d^{1/3}}{16 (\log d)^{1/3}}|X|.
$$
This completes the proof of the corollary.
\vspace{0.3cm}

\noindent
By the last Corollary and Hall's theorem, one can assign a set of 
$\frac{1}{16}d^{1/3}/(\log d)^{1/3}$ members of 
$L \cup (F \setminus T)$ to any
element of  $(C \setminus H) \cup T$, so that no 
member of $L \cup (F \setminus T)$ 
is assigned more than once. This makes all elements of 
$(C \setminus H) \cup T$ centers of
large vertex disjoint stars. Adding to these stars 
the stars whose centers
are the elements of $H$ and whose leaves are those of $S$,
we may apply Claim \ref{c32} 
to conclude that in case there are any unassigned vertices left in
$L$ we can connect each of them to one of the 
existing centers.
Similarly, (\ref{e31})  
and the definition of rule (a)  
do the same job for unassigned vertices in $F \setminus T$.
Thus, we get a star factor in which each star has at least
$\frac{1}{25}d^{1/3}/(\log d)^{1/3}$ leaves. This completes the proof.
\qed

\section{Concluding remarks and open problems}

We have shown that for every positive integer $g$ there is 
an integer $d$ so that any graph with minimum degree at least $d$
contains a star factor in which every component has at least
$g$ edges. Let $d(g)$ denote the minimum number $d$ for which this 
holds. Our main result shows that
$d(g) \leq O( g^3 \log g)$, while the construction described in the
proof of Theorem \ref{t12} implies that $d(g) \geq \Omega( g \log g)$.
It seems plausible to conjecture that $d(g)=\Theta( g \log g)$, but 
this remains open. It will be interesting to determine
$d(g)$ precisely (or estimate it more accurately) 
for small  values of $g$, like $g=2$ or $3$. Our 
proof, even if we try to optimize the constants in it, will
yield only some crude upper bounds that are certainly far from 
being tight. It is worth noting, however, that even for showing
that $d(2)$ is finite, we do not know
any proof simpler than the one given here for the general case. 
On the other hand, a random 3-regular graph contains a Hamilton 
cycle with probability tending to 1. Hence, if it has number of 
vertices divisible by 3, it contains a spanning factor of stars of 
two edges each. Similarly, a random 4-regular graph  having 
number of vertices divisible by 4 contains a spanning factor of 
stars of three edges each  with probability tending to 
1~\cite{AW}. Immediately from  contiguity results discussed 
in~\cite{W}, the same statements are true if we change 4-regular 
to $d$-regular for any $d\ge 4$.

Our proof, together with
the algorithmic version of the local lemma proved by Beck
in \cite{Be} (see also \cite{Al1}), and any efficient algorithm
for bipartite matching, show that the proof here can be converted to 
a deterministic, polynomial time algorithm that finds, in
any given input graph with minimum degree at least $d$, a star
factor in which every star is of size at least
$\Omega((d/ \log d)^{1/3})$. We omit the details.

There are several known results that show that any graph with
large minimum degree contains a spanning tree with many leaves,
see \cite{KW}, \cite{GW}, \cite{CWY}. In particular, it is known
(and not difficult)
that any graph with minimum degree $d$ and $n$ vertices contains
a spanning tree with at least $n-O\big( n (\log d)/d\big)$ leaves.
A related question to the one considered here is whether it is true
that any graph with large minimum degree contains a spanning tree in
which all non-leaf vertices have large degrees. Specifically, 
is there an absolute positive constant $c$ so that
any graph with minimum degree at least $d$ contains
a spanning tree in which the degree of any non-leaf is at least
$c d/\log d$ ? Another intriguing question is the following possible
extension of the main result here. Is it true that the edges of
any graph $G$ with minimum degree $d$ can be partitioned into
pairwise disjoint sets, so that each set forms a spanning forest
of $G$ in which every component is of size at least $h(d)$,
where $h(d)$ tends to infinity  with $d$ ? A related result
is proved in \cite{AMR}, but  the proof of the  last statement,
if true, seems to require additional ideas.


\end{document}